\documentclass[12pt,a4paper]{article}
\usepackage{amsmath,amsthm,amsfonts,latexsym,amscd,amssymb,enumerate}

\usepackage[hypertex,backref]{hyperref}

\newtheorem{theorem}{Theorem}[section]
\newtheorem{thm}[theorem]{Theorem}
\newtheorem{prop}[theorem]{Proposition}
\newtheorem{lem}[theorem]{Lemma}
\newtheorem{cor}[theorem]{Corollary}

\theoremstyle{remark}
 
\newtheorem{rem}[theorem]{Remark}

\theoremstyle{definition}

\newtheorem{corollary}[theorem]{Corollary}

\theoremstyle{definition}
\newtheorem{definition}[theorem]{Definition}

\theoremstyle{remark}
\newtheorem{remark}[theorem]{Remark}

\newcommand{\innerprod}[1]{\langle #1 \rangle}
\DeclareMathOperator{\ind}{ind}

\newcommand{\complexs}{\mathbb{C}}
\newcommand{\naturals}{\mathbb{N}}
\newcommand{\integers}{\mathbb{Z}}
\newcommand{\rationals}{\mathbb{Q}}

\newcommand{\tensor}{\otimes}
\newcommand{\iso}{\cong}
\newcommand{\norm}[1]{\left\lVert#1\right\rVert}

\DeclareMathOperator{\End}{End}
\DeclareMathOperator{\ch}{{\rm ch}}

\DeclareMathOperator{\U}{U}
\DeclareMathOperator{\Spc}{Spin^c}

\DeclareMathOperator{\Mat}{Mat}

\DeclareMathOperator{\vol}{vol}
\DeclareMathOperator{\Todd}{\mathcal{T}}
\newcommand{\K}{\mathbb{ K}}
\newcommand{\C}{\mathbb{ C}}
\newcommand{\B}{\mathbb{ B}}

\newcommand{\N}{\mathbb{ N}}

\newcommand{\Z}{\mathbb{ Z}}

\newcommand{\Hom}{\operatorname{Hom}}

\newcommand{\Q}{\mathbb{ Q}}
\newcommand{\R}{\mathbb{ R}}

\newcommand{\ra}{{\rightarrow}}

\newcommand{\id}{\operatorname{id}}

{\catcode`@=11\global\let\c@equation=\c@theorem}

\begin{document}

\title{Enlargeability and index theory}

\author{B.~Hanke and T.~Schick}

\maketitle

\begin{abstract} Let $M$ be a closed enlargeable spin manifold.
We show non-triviality of the universal index obstruction in the
$K$-theory of the maximal $C^*$-algebra of the fundamental group of
$M$. Our proof is independent from the injectivity of the Baum-Connes assembly
map for $\pi_1(M)$ and relies on the construction of a certain infinite
dimensional flat vector bundle out of a sequence of finite
dimensional vector bundles on $M$ whose curvatures tend to zero.

Besides the well known fact that $M$ does
not carry a metric with positive scalar curvature, our results
imply that the classifying map $M \to B \pi_1(M)$
sends the fundamental class of $M$ to a nontrivial homology class in
$H_n(B \pi_1(M) ; \Q)$. This answers a question of Burghelea (1983).
\end{abstract}

\section{Introduction}

\subsection{Enlargeability and the universal index obstruction}
\label{sec:enlarg-univ-index}

For a closed spin manifold $M^n$, Rosenberg in \cite{Ro3} constructs an index
\[
    \alpha^{\R}_{max}(M) \in KO_n(C^*_{\rm max,\R}\pi_1(M))
\]
in the $K$-theory of the (maximal) real $C^*$-algebra of the fundamental
group of $M$. By the Lichnerowicz-Schr{\"o}dinger-Weitzenb{\"o}ck formula this index is zero if $M$ admits
a metric of positive scalar curvature. The Gromov-Lawson-Rosenberg
conjecture states that, conversely, the vanishing of $\alpha(M)$
implies that $M$ admits such a metric, if $n \geq 5$. By a
result of one of the authors, this conjecture
is known to be false in general \cite{Sch0}. But a stable version of
this conjecture is  true, if the Baum-Connes assembly map
\[
   \mu: KO^{\pi_1(M)}_*(\underline{E}\pi_1(M)) \to KO_*(C^*_{max,\R} \pi_1(M))
\]
is injective \cite{St}. The proof of this (and related results) is based on the
existence of a natural map $D:KO_*(M) \to KO_*^{\pi_1(M)}(\underline{E} \pi_1(M))$
into the equivariant $K$-homology of the classifying space for
$\pi_1(M)$-actions with finite isotropy and of a factorization
\[
   KO_n(M) \stackrel{D}{\to} KO^{\pi_1(M)}_n(\underline{E}\pi_1(M))
\stackrel{\mu}{\longrightarrow}
   KO_n(C^*_{max,\R} \pi_1(M))
\]
which sends the $KO$-fundamental class $[M] \in KO_n(M)$ to 
$\alpha^{\R}_{max}(M)$.
Therefore, if $\alpha^{\R}_{max}(M) = 0$ and $\mu$ is injective, one knows
that $D([M])=0$ and this situation can be analyzed by algebraic
topological means. (Actually, Stephan Stolz is using the reduced group
$C^*$-algebra, compare the discussion in Section \ref{sec:different-c-indices}.)

In this paper, we describe a new method to detect non-vanishing of
this universal index obstruction in a nontrivial case. This is
independent of the injectivity of the Baum-Connes map. For convenience,
we study the complex K-theory index element $    \alpha_{max}(M) $
in the $K$-theory of the maximal complex $C^*$-algebra of $\pi_1(M)$.
The usual Lichnerowicz argument shows that $\alpha_{max}(M) = 0$ if $M$ admits
a metric of positive scalar curvature. 

In the first part of our paper, we prove a weak  converse to this statement.
Recall:

\begin{definition} A closed oriented manifold $M^n$ is called
  {\em enlargeable} if the following holds:
  Fix some Riemannian metric $g$ on $M$. Then, for all $\epsilon > 0$,
  there is a finite cover $\overline{M}$ of $M$  and an
  $\epsilon$-contracting map $(\overline{M}, \overline{g}) \to (S^n,
  g_0)$ of non-zero degree, where $\overline{g}$ is induced by $g$ and
  $g_0$ is the standard metric on $S^n$.

  $M$ is called \emph{area-enlargeable} if in the above
  $\epsilon$-contracting is replaced by \emph{$\epsilon$-area
    contracting}. Here, a map $f\colon M\to N$ between two
  $n$-dimensional Riemannian manifolds is called $\epsilon$-area
  contracting if $\norm{\Lambda^2 T_xf}\le \epsilon$ for each $x\in
  M$, where $\Lambda^2T_xf$ is the induced map on the second exterior
  power $\Lambda^2T_x M\to \Lambda^2 T_{f(x)}N$ with norm induced by
    the Riemannian metrics.
\end{definition}

Note that every enlargeable manifolds is area-enlargeable, but that
the converse might not be true.

We remark that in contrast to the definition
in \cite{GL}, we do not require that the
covers $\overline{M}$ necessarily admit spin structures.

\begin{thm} \label{mainth} Let $M$ be an enlargeable or
  area-enlargeable spin manifold. Then
\[
   \alpha_{max}(M) \neq 0 \in K_n(C^*_{max} \pi_1(M)) \, .
\]
\end{thm}

By a result of Gromov and Lawson \cite{GL}, enlargeable spin manifolds do
not admit metrics of positive scalar curvature. Recall the question posed in the second
paragraph of the introduction to the article \cite{Ro1}: ``Nevertheless, it is
not clear, if their results always imply ours or vice versa.''

Our paper gives a complete answer in one direction: 
if $M$ is spin,
the index obstruction $\alpha_{max}(M)$ completely subsumes the
enlargeability (and area-enlargeability) obstruction to positive scalar curvature of Gromov and Lawson.

For the applications to positive scalar curvature, we restrict our
discussion to spin manifolds $M$ in order
to keep the exposition transparent. It should not be too hard to
extend the methods and results to the case where only the universal
cover of $M$ admits a spin structure. 

\subsection{Flat bundles of $C^*$-modules}
\label{sec:flat-bundles-c}

The idea of our proof can be summarized as follows. We construct
a $C^*$-algebra morphism
\begin{equation*}
    \phi: C^*_{max}\pi_1(M)  \to Q\label{eq:define_phi}
\end{equation*}
where $Q$ is a (complex) $C^*$-algebra whose $K$-theory can be explicitely calculated,
and then we study the image of $\alpha_{max}(M)$ under the induced
map in $K$-homology.

The map $\phi$ results from the holonomy representation of $\pi_1(M)$
associated to an infinite dimensional flat bundle on $M$ which is
obtained in the following way: Because $M$ is enlargeable or area-enlargeable, there is a sequence
$E_i \to M$ of (finite dimensional) unitary vector bundles
with connections whose curvatures tend to zero, but whose Chern characters are
nontrivial. We then construct an infinite dimensional smooth bundle $V \to M$
with connection and with the following property: The fiber over $p \in M$ consists of bounded sequences
$(v_1, v_2 , \ldots)$ with $v_i \in (E_i)_p$ and the connection restricts
to the given connection of $E_i$ on each ``block''. We denote by $W \subset V$
the subbundle consisting of sequences tending to zero. The
$\End(V)$-valued curvature form on $V$
sends $V$ to $W$ by the asymptotic curvature property of the sequence
$(E_i)$. Hence the quotient bundle
$V/W \to M$ with the induced connection is flat.

However, this bundle still encodes the asymptotic non-triviality
of the Chern characters of the original bundles and hence
the index of the Dirac operator on $M$
twisted with this bundle is nontrivial. This index can be expressed in terms
of the collection of indices of the Dirac operator twisted with $E_i$, $i \in \N$. 

It should be noted that the precise argument in Section
\ref{sec:assembling} needed to construct the bundle $V$
requires a considerable amount of care.

In this respect, we realize the idea formulated at the end of the introduction
to \cite{GLIHES}: ``Passing to the limit, one might expect to find an
interesting infinite dimensional, {\em flat} bundle $E_0$ over the original
manifold, so that one could apply the Bochner method directly to the Dirac
operator with coefficients in $E_0$''. In our case, the role of $E_0$ is
played by the bundle $V/W \to M$.

\subsection{Almost flat bundles and almost representations}
\label{sec:almost-flat-bundles}

Our construction can be seen in relation to the notions
of almost flat bundles as studied by Connes-Gromov-Moscovici \cite{CGM}
and of almost representations as studied by Mishchenko and his coauthors
(compare e.g.~\cite{MM}). Heath Emerson informed us that he and Jerry
Kaminker plan to carry out a systematic study of these notions in the context
of the Baum-Connes conjecture. Contrary to the definitions used
in the mentioned sources we do not require the different bundles in
the almost flat sequence to define the same $K$-theory class
or the ``not quite representations'' induced by such a sequence
to be related in any way. Keeping this flexibility throughout
the argument enables us to  prove the general statement
of Theorem \ref{mainth}.

\subsection{The different $C^*$-indices}
\label{sec:different-c-indices}

In our paper, we use complex $C^*$-algebras, because
this avoids some technicalities and is sufficient for the applications
that we have in mind. The respective real versions of our theorems can
be shown similarly. 

Here we want to compare the reduced to the
maximal index, and the real to the complex version. In recent
literature on the positive scalar curvature question, in most cases
the real reduced index $\alpha_{red}^{\R}(M)\in
KO_n(C^*_{red,\R}\pi_1(M))$ is used, whereas Rosenberg \cite{Ro1,Ro3}
uses $\alpha_{max}(M)$ and $\alpha_{max}^{\R}(M)$.

 Note
that for any discrete group $\pi$, we have canonical maps
\[
\begin{CD}
  C^*_{max,\R} \pi @>{\omega^{\R} }>> C^*_{red,\R} \pi\\
  @VVV @VVV\\
  C^*_{max} \pi @>{\omega }>> C^*_{red} \pi\\
\end{CD}
\]
and a commutative diagram
\begin{equation}
\begin{CD}
KO_*^{\pi}(\underline{E}\pi) @>{\mu_{max}^{\R}}>> KO_*(C^*_{max,\R} \pi)
@>{\omega^\R}>> KO_*(C^*_{red,\R} \pi) \\ 
@VVV      @VVV @VVV \\
K_*^{\pi}(\underline{E} \pi) @>{\mu_{max}}>> K_*(C^*_{max}\pi)
@>{\omega}>> K_*(C^*_{red}\pi) \, 
\end{CD}\label{eq:cstardiagram}
\end{equation}
where the vertical maps are given by complexification and the
horizontal compositions are the reduced analytic assembly maps $\mu$
or $\mu^\R$.
If $\pi = \pi_1(M)$, the map $\omega^{?}_*$ sends $\alpha_{max}^{?}(M)$ to
$\alpha_{red}^{?}(M)$ (this is true for the real and the complex 
version), and the complexification maps send
$\alpha_{?}^{\R}(M)$ to $\alpha_{?}(M)$ (for the max and red
version). In particular, vanishing of $\alpha_{max}^{\R}(M)$ implies vanishing
of all the other index invariants. Consequently, following Rosenberg,
one should formulate the
Gromov-Lawson-Rosenberg conjecture using $\alpha_{max}^{\R}$. We
point out that this conjecture holds
stably by the result of Stephan Stolz cited above, if the Baum-Connes map $\mu_{max}^\R$ is
injective. It  could possibly happen that
$\alpha_{max}^{\R}(M) \neq 0$ whereas $\alpha_{red}^{\R}(M) = 0$. However,
this would imply that a number of important conjectures are
wrong, most notably that the (reduced) Baum-Connes assembly map is not always
injective.

The maximal $C^*$-algebra has much better functorial properties than
the reduced one, a fact that we are using in the construction of the
homomorphism $\phi$ alluded to in \eqref{eq:define_phi}. 

In view of these considerations and the results presented in this
paper, we propose to always use the obstruction $\alpha^{\R}_{max}(M)$
to the existence of positive scalar curvature metrics instead of
$\alpha^{\R}_{red}(M)$. Note that the two obstructions
are equivalent, if the Baum-Connes map is injective. Note also that
the two obstructions coincide if the fundamental group $\pi$ is K-amenable
(in particular if it is amenable) because in this case the map
$K_*(C^*_{max}\pi)\to K_*(C^*_{red}\pi)$ is an isomorphism.

In this paper, we show non-vanishing of the complex version
$\alpha_{max}(M)$ under an enlargeability assumption, which implies by
\eqref{eq:cstardiagram} non-vanishing of $\alpha^{\R}_{max}(M)$.

\subsection{Enlargeability and the  fundamental class}
\label{sec:enlarg-homol-fund}

Turning to another application of our methods, we show

\begin{thm} \label{burg} Let $M$ be a closed enlargeable or area-enlargeable manifold,
$f:M \to B \pi_1(M)$ classify the universal cover of $M$
and $[M] \in H_n(M ; \Q)$ be the fundamental class. Then
\[
   f_{*}([M]) \neq 0 \in H_n(B \pi_1(M) ; \Q)
\]
\end{thm}

This theorem implies an affirmative answer to a question of
Burghelea \cite[Problem 11.1]{Schu}. We emphasize that
contrary to the original formulation of Burghelea's question, no spin
assumption on $M$ or its universal cover is required.
It is somewhat remarkable that
we prove Theorem \ref{burg} by making a detour through the $K$-theory
of $C^*$-algebras and an assembly map.

For (lenth)-enlargeable manifolds, one can use coarse geometry methods to
get a shorter proof of this result. It would be interesting to extend these coarse methods to
the area-enlargeable case.

For $n \geq 4$, we will construct closed oriented manifolds 
$M^n$ (that may be chosen to be spin) whose 
fundamental classes are sent to nontrivial classes in $H_n(B\pi_1(M);\Q)$, but 
which are not area-enlargeable. In this respect, a converse of Theorem 
\ref{burg} does not hold to be true.

Because the proof of Theorem \ref{burg} is based on an
analysis of general Dirac type operators on $M$, the necessary index
theory in Section \ref{sec:indextheory} will be developed
in the required generality.

\begin{corollary}
  It follows from Theorem \ref{burg} that every (area)-enlargeable manifold $M$ is essential in the sense
  of Gromov's \cite{Gromov} and therefore its $1$-systole satisfies
  Gromov's main inequality 
  \begin{equation*}
    sys_1(M) \le c(n) \vol(M)^{1/n}.
  \end{equation*}
  In particular, such an $M$ has a non-contractible closed geodesic of
  length at most $c(n) \vol(M)^{1/n}$. Here, $n=\dim(M)$ and
  $c(n)>0$ is a constant which depends only on this dimension.
\end{corollary}

\section{Assembling almost flat bundles}
\label{sec:assembling}

Let $M$ be a closed smooth $n$-dimensional Riemannian manifold,
let $d_i$, $1 \leq i < \infty$, be a sequence of natural numbers
and let $(P_i, \nabla_i)_{i \in \N }$ be an {\em almost flat} sequence
of principal $\U(d_i)$ bundles over $M$ equipped with $\U(d_i)$-connections
$\nabla_i$. By definition, this means that the curvature $2$-forms
\[
    \Omega_i \in \Omega^2(M;\mathfrak{u}(d_i)  )
\]
associated to $\nabla_i$ vanish asymptotically with respect to
the maximum norm on the unit sphere bundle in $\Lambda^2 M$ and the operator norm
on each $\mathfrak{u}(d_i) \subset \Mat(d_i) := \C^{d_i \times d_i}$, i.e.
\[
    \lim_{i \to \infty} \| \Omega_i \| = 0 \, .
\]
Let $\K$ denote the $C^*$-algebra of compact operators on $l^2(\naturals)$. We
choose embeddings of complex $C^*$-algebras
\[
  \gamma_i\colon \Mat(\C, d_i) \hookrightarrow \K \, .
\]
Hence, each $P_i$ has an  associated bundle
\[
    F_i := P_i \times_{U(d_i)} \K
\]
consisting of projective right $\K$-modules with one generator and which is equipped with a $\K$-linear connection
\[
   \nabla_i:\Gamma(F_i)  \to  \Gamma(T^*M \otimes F_i ).
\]
Note that, since the map $U(d_i)\to \K$ is not unital, the fibers of $F_i$ are
\emph{not} free, but are isomorphic as right $\K$-modules to $q_i\K$,
where $q_i=\gamma_i(1)$.
Because the structure group of $F_i$ is $\U(d_i)$, each bundle $F_i$ has the
structure of a $\K$-Hilbert bundle induced by the inner product
\[
   \K \times \K \to \K \, , (X,Y) \mapsto X^* Y
\]
on each fiber. The connections $\nabla_i$ are compatible with these
inner products. Let $A$ be the complex unital $C^*$-algebra of norm bounded sequences
\[
     (a_i)_{i \in \N} \in \prod_{i=1}^{\infty} \K  \, .
\]
For $i \in \N$, we denote by $A_i \subset A$ the subalgebra of sequences such
that all but the $i$th entry vanish. The algebra $A_i$ can
be identified with $\K$. Define the element (being a projection)
\begin{equation*}
  q:=(q_i)_{i\in\N} \in A;\qquad q_i=\gamma_i(1).
\end{equation*}
The following theorem says that the bundles $F_i$ can be assembled to a
smooth bundle of right Hilbert $A$-modules in a particularly nice way.
For the necessary background concerning Hilbert module bundles, we
refer to \cite{Sch}.

\begin{thm} \label{main} There is a smooth Hilbert $A$-module bundle $V \to M$,
each fiber of $V$ being a finitely generated projective right $A$-module, together with an
$A$-linear metric connection
\[
   \nabla^V :\Gamma(V) \to \Gamma(T^*M \otimes V)
\]
such that the following holds:
\begin{itemize}
\item  For $i \in \N$, let $V_{i}$ be the
subbundle $V \cdot A_{i} \subset V$. Then $V_{i}$ (``the $i$th
block of $V$'') is isomorphic to $F_{i}$  (as a $\K$-Hilbert
bundle).
\item The connection $\nabla$ preserves the subbundles $V_i$.
\item Let $\nabla^V_i$ be the connection induced on $V_i$ by $\nabla^V$
and let $\Omega^V_i$ be the corresponding curvature form in
$\Gamma(\Lambda^2 M \otimes \End_{A_i}(V_i))$. Then
\[
    \lim_{i \to \infty} \| \Omega^V_i \| = 0 \, .
\]
\end{itemize}
\end{thm}

The remainder of this section is devoted to the construction of $V$. At
first, we obtain a Lipschitz-Hilbert-$A$-module
bundle $L \to M$ that will be approximated by a smooth bundle
$V \to M$. After this has been done, the bundle $V$ will
be equipped with a connection $\nabla^V$ as stated in Theorem \ref{main}.

\medskip

Let
\[
   D^n : =  \{ (x_1, \ldots, x_n) \in \R^n ~|~ 0 \leq x_i \leq 1 \} \subset
\R^n
\]
be the standard $n$-dimensional cube and let $(\phi_j)_{j \in J}$ be a
finite family of diffeomorphisms \footnote{in the
sense that each $\phi_j$ extends to a diffeomorphism from
an open neighborhood of $W_j \subset M$ to an open neighborhood of
$D^{n} \subset \R^n$}
\[
    M \supset W_j \stackrel{\phi_j}{\to} D^n
\]
such that
\[
    M \subset \bigcup_{j \in J} \stackrel{\circ}{W}_j \, .
    \]
Identify each of the $W_j$ with $D^n$, using $\phi_j$.
In order to obtain the bundle $L$, we construct trivializations
\[
   \psi_{i,j}: F_i|_{W_j} \cong D^n \times q_i\K
\]
for all $i\in \N$ and $j \in J$ as follows: Choose a $\K$-linear isomorphism
\[
    \psi_{i,j}: F_i|_{(0,0, \ldots, 0)} \cong q_i\K .
\]
The map $\psi_{i,j}$ can be extended to a unique isomorphism
of smooth $\K$-module bundles
\[
    \psi_{i,j}:F_i|_{[0,1] \times 0 \times \ldots \times 0} \cong
  ([0,1] \times 0 \times \ldots \times 0) \times q_i\K
\]
such that the constant sections
\[
  [0,1] \times 0 \times  \ldots \times 0 \to ([0,1] \times 0 \times \ldots
\times 0) \times q_i\K
\]
are parallel with respect to $\nabla$. Inductively,
we assume that $\psi_{i,j}$ has already been defined on
\[
   F_i|_{D^k \times 0 \times \ldots \times 0} \, .
\]
Then $\psi_{i,j}$ can be extended to a unique isomorphism of
smooth $\K$-module bundles
\[
    F_i|_{D^{k+1} \times 0 \times \ldots \times 0} \cong D^{k+1} \times
    q_i\K
\]
such that the covariant derivative along  the tangent vector field
\[
    \frac{\partial}{\partial x_{k+1}} \in \Gamma(TD^{k+1})
\]
of each constant section $D^{k+1} \to D^{k+1} \times q_i\K$ vanishes.

\begin{definition}
We denote by
\[
   \omega_{i,j} \in \Gamma(T^* D^n \otimes \End_{\K}(q_i\K)) \iso
   \Omega^1(D^n; q_i\K q_i)
\]
the connection $1$-form induced on $D^n \times q_i\K$ by $\nabla_i$ and
$\psi_j$. Note that the right $\K$-module endomorphisms of
$q_i\K$ are canonically isomorphic to the unital $C^*$-algebra $q_i\K
q_i$.

Furthermore, we denote by $\| \omega_{i,j}
\|$ the $L^{\infty}$-norm of $\omega_{i,j}$ induced by the usual Euclidean
metric on $D^n$ and the operator norm on $\End_{\K}(q_i \K)$.
\end{definition}

We will show now that the special construction of the trivializations
$\psi_{i,j}$ ensures
that we have upper bounds for $\| \omega_{i,j} \|$. Let
\[
   \eta_{i,j} = d \omega_{i,j} - \omega_{i,j} \wedge \omega_{i,j} \in
\Gamma(\Lambda^2 D^n \otimes  \End_{\K} ( D^n \times q_i\K)) =
\Omega^2(D^n;q_i\K q_i)
\]
be the curvature $2$-form on $D^n$ induced by $\psi_{i,j}$ and $\nabla_i$.

\begin{lem} \label{estimate} For each $i$ and $j$, we have
\[
      \| \omega_{i,j} \|  \leq n \cdot  \| \eta_{i,j} \| \, ,
      \]
      where $n=\dim(M)$.
\end{lem}

\begin{proof} For brevity, we drop the indices $i$ and $j$  and
abbreviate $\frac{\partial}{\partial x_{\nu}}$ by $\partial_{\nu}$.
By construction of the trivialization $\psi$, we have
\[
 \omega_{(x_1, \ldots, x_k, 0 \ldots, 0)}(\partial_{\nu}) = 0 \, ,
\]
if $\nu \geq k$. Now, if $k > \nu$, we get
\begin{eqnarray*}
 \lefteqn{\| \omega_{(x_1,  \ldots, x_k, 0, \ldots, 0)}(\partial_{\nu}) \| = } \\
 &  &  = \| \omega_{(x_1, \ldots, x_{k-1}, 0, \ldots, 0)}(\partial_{\nu}) +    \int_0^{x_k}   d\omega_{(x_1, \ldots, x_{k-1}, t, 0, \ldots, 0)}
 (\partial_k, \partial_{\nu}) \| \\
 &  & \le \| \omega_{(x_1, \ldots, x_{k-1}, 0, \ldots, 0)}(\partial_{\nu}) \| + \int_0^{x_k} \| \eta_{(x_1, \ldots, t, 0, \ldots, 0)}(\partial_k, \partial_{\nu}) \| \\
& & \leq \| \omega_{(x_1, \ldots, x_{k-1}, 0, \ldots, 0)}(\partial_{\nu})\|  + \|\eta_{i,j} \| \cdot |x_k|.
\end{eqnarray*}
The second inequality uses the fact that
\[
(\omega \wedge \omega)_{(x_1, \ldots, x_{k-1}, t, 0, \ldots, 0)}
(\partial_k, \partial_{\nu}) = 0
\]
by construction of the trivialization $\psi_{i,j}$. Because
$\omega_{(x_1, \ldots, x_{\nu}, 0 , \ldots, 0)}(\partial_{\nu}) = 0$,
we see inductively that
\[
 \| \omega_{(x_1, \ldots, x_k, 0, \ldots, 0)}\| \leq \| \eta_{i,j} \| \cdot (|x_{\nu+1}| + \ldots + |x_k|)
\leq n  \cdot \| \eta_{i,j} \| \, .
\]
\end{proof}

\begin{remark}
  By assumption, the bundles $P_i$ form an almost flat sequence of
  bundles. Consequently, the supremum norms $\norm{\eta_{i,j}}$ (which
  we know do not depend on the particular trivializations) have an
  upper bound, and by Lemma \ref{estimate} the same is true for the
  $\norm{\omega_{i,j}}$.
\end{remark}

\begin{lem} \label{easyestimate} Let $l \geq 0$. Then there is
a constant $C(l)$ (independent of $i,j$) such that
if
\[
     \phi:[0,1] \to D^n \times q_i\K
\]
is a parallel vector field (with respect to $\omega_j$) along a
piecewise smooth path $\gamma: [0,1] \to D^n$ of length $l(\gamma) \leq l$
(measured with respect to the usual metric on $D^n$), then
\[
 \| \phi(1) - \phi(0) \| \leq C(l) \cdot \| \omega_{i,j} \| \cdot
 l(\gamma)\cdot \norm{\phi(0)}
\]
for all $i,j$. The constant $C(l)$ depends on the supremum of all the $\norm{\omega_{i,j}}$.
\end{lem}

\begin{proof} Since the bundle $D^n\times q_i\K$ is trivial, we
  consider the section $\phi$ as a path $[0,1] \to q_i\K$. It satisfies the
differential equation
\[
 \phi'(t) + \left( (\omega_{i,j})_{\gamma(t)}(\dot{\gamma}(t))\right)\cdot
 \phi(t) = 0  
 \]
and  
it follows that
\[
  \norm{\phi(1)-\phi(0)}\le
  \exp\left(2l(\gamma)\norm{\omega_{i,j}}\right)\cdot \norm{\phi(0)} \, .
\]
The function $\exp:q_i\K q_i \to q_i \K q_i$ is uniformly Lipschitz continuous
on each bounded neighborhood of $0$. Hence, the proof is complete.
\end{proof}

These estimates allow for the following important implication.
For $\alpha, \beta \in J$, $i\in\N$, we denote by
  \begin{equation*}
   \phi_{\alpha,\beta} : \phi_{\alpha}(W_{\alpha} \cap W_{\beta})  \to 
                         \phi_{\beta}(W_{\alpha} \cap W_{\beta}) 
                       \end{equation*}
the transition function for the charts $\phi_i$ of our manifold $M$, and
\begin{equation*}
\psi_{\alpha,\beta,i} : \psi_{\alpha} (F_i|_{W_{\alpha} \cap
     W_{\beta}})  \to \psi_{\beta} (F_i|_{W_{\alpha} \cap
     W_{\beta}})
 \end{equation*}
the transition function for the trivializations of the bundles $F_i$.

\begin{prop} \label{Lipschitz} There is a constant $C \in \R$ such that (independent
of the particular smooth bundle in the almost flat sequence) the following holds
for all $\alpha$ and $\beta$: Considering $\psi_{\alpha, \beta,i}$
(i.e.~the (smooth) transition function for the $i$-th vector bundle)
as a function
\[
   \psi: W_{\alpha} \cap W_{\beta} \to q_i\K q_i \, ,
\]
we have $\|D \psi(x)\| \leq C$ for all  $x \in W_{\alpha} \cap W_{\beta}$.
\end{prop}

\begin{proof} Let
\[
 x = (x_1, \ldots, x_n) \in \phi_{\alpha} (W_{\alpha} \cap W_{\beta}) \cap \stackrel{\circ}{D^n}
\]
and let $1 \leq \nu \leq n$.   We have to study the function
\[
   f: (-\epsilon, \epsilon) \to q_i\K q_i
\]
which is the restriction of $\psi$ to a certain path and
defined by the property that
\[
   \psi_{\beta}\psi_{\alpha}^{-1} (x + t e_{\nu}, v) =
   (\phi_{\alpha, \beta}(x+ te_{\nu}),  f(t) \cdot v)
\]
for all $v \in q_i\K$ and all $t \in (-\epsilon, \epsilon)$
(where $\epsilon$ is  sufficiently small). For all $t \neq 0$, the
fact that parallel transport is isometric implies that each $f(t)$ is
unitary and hence
\[
   \| \frac{f(t) - f(0)}{t} \|  = \| \frac{f(t) \circ f(0)^{-1} - \id}{t}\| \, .
\]
For $v \in q_i\K$, the element $(f(t) \circ f(0)^{-1}) (v)\in q_i\K$
can be  constructed as follows:  Consider the path
\[
   \gamma: [0, t] \to D^n \, , ~ \xi \mapsto x + \xi e_{\nu} \, .
\]
Now parallel transport the element $v$ along $\gamma^{-1}$ using the connection
$\omega_{\alpha}$ to get $w \in q_i\K$ and then parallel transport this element
$w$ along $\phi_{\alpha, \beta} \circ \gamma$ using the connection
$\omega_{\beta}$. This works since, in terms of the bundle $F_i$, this means we use the
same parallel transport (given by $\omega_\alpha$ and $\omega_\beta$
in the two trivializations) to transport a given vector in the fiber
of $\gamma(t)$ to the fiber of $\gamma(0)$, where, in terms of the
trivialization, they are identified using $f(0)$.

By Lemmas \ref{easyestimate} and \ref{estimate}, together
with the fact that the curvatures of the bundles $F_i$
are universally bounded, the norms
\[
   \| v - w \| {\rm~and~}  \| w - f(t) \circ f(0)^{-1}(v) \|
\]
are bounded up to a universal constant by the length of $\gamma$ or
$\phi_{\alpha, \beta} \circ \gamma$, respectively, and hence are bounded by $C' \cdot t$,
where $C'$ is a constant independent of $i$, $\alpha$ and $\beta$ (note that
$\|D\phi_{\alpha,\beta}\|$ is uniformly bounded as the supremum of
finitely many compactly supported continuous functions).
This implies the assertion of the proposition with $C := 2C'$.
\end{proof}

We call a continuous Banach-space bundle
\[
   F \hookrightarrow L \to M
\]
(where the typical fiber $F$ is a complex Banach space) a {\em Lipschitz bundle}
if the following holds: There is an open covering  $(U_j)_{j \in J}$ of $M$
and there are trivializations
\[
      L|_{U_j} \cong U_j \times F
\]
so that the associated transition functions
\[
  U_{\alpha} \cap U_{\beta} \to \End(F)
\]
are locally Lipschitz continuous,  $\End(F,F)$ being equipped with
the operator norm

In this context, Proposition \ref{Lipschitz} can be summarized as follows:

\begin{thm} \label{L} The bundles $F_i \to M$ can be assembled to a Lipschitz
Hilbert $A$-module bundle
\[
   L \to M \,
   \]
   with typical fiber $qA$.
The bundles $L \cdot A_i$ all have a smooth structure
compatible with the induced Lipschitz structure and are  isomorphic
to $F_i$.
\end{thm}

In the following, we will use results about Hilbert
$A$-module bundles as explained in
\cite{Sch} where the role of smooth structures of Hilbert
$A$-module bundles is quite carefully explained. Here, we will
frequently use Lipschitz structures of such bundles (i.e.~the
transition functions of a Lipschitz atlas are (locally) Lipschitz
continuous). It is straightforward to check that all the results
described in
\cite{Sch} we are using here carry over immediately to the Lipschitz
category. 

In order to construct the bundle $V$ described in Theorem \ref{main}, 
we use \cite[Theorem 3.14]{Sch}  and write the bundle
$L \to M$ as a subbundle of a trivial Hilbert $A$-module bundle
\[
   M \times A^k \to M \, .
\]
Hence, $L$ is the image of a projection valued
\emph{Lipschitz continuous} section $\phi$ of this bundle. The section $\phi$ can
be approximated arbitrarily close (in the operator norm of $\Hom_A(A^k, A^k)
= A^{k \times k}$ and the maximum norm on $M$) by a smooth
projection valued section.
The resulting bundle $V$ (consisting of the images of these projections) is a smooth Hilbert  $A$-module bundle. We choose the approximation close
enough such that $V$ is Lipschitz isomorphic as a Hilbert $A$-module
bundle to the bundle $L$, in particular it also has  typical fiber
isomorphic to $qA$ (cf. \cite[Lemma 3.12.]{Sch}).

By the algebraic structure of $A$, also $V$ has ``blocks'' $V
\cdot A_i$ and, being an $A$-module bundle isomorphism, the
isomorphism between $V$ and $L$ maps the blocks $V
\cdot A_i$ to $L\cdot A_i$. By construction of $V$, this way we obtain
\emph{smooth} Hilbert
$A_i$-module bundle isomorphisms
\[
    V \cdot A_i \cong F_i \, .
\]
The
trivializations $\psi_{i,j}$ assemble to a  Lipschitz continuous trivialization
\[
   L|_{W_j} \cong W_j \times qA \, .
   \]
On the other hand, we can choose smooth Hilbert $A$-module
bundle trivializations
\begin{equation*}
   V|_{W_j} \cong W_j \times qA;
 \end{equation*}
because the $W_j$ are contractible and the typical
fiber of $V$ is isomorphic to $qA$.

Observe that $\End_A(qA) \cong qAq$, where $\End_A(qA)$ denotes the right
$A$-module endomorphisms, and $qAq$ acts by left
multiplication. The isomorphism
$L \cong V$ can hence be described by a Lipschitz continuous map
\[
  \tau_j :  D^n \to qAq \,
\]
with values in the unitary group of $qAq$.

The desired connection on $V$ is now constructed as follows.
Let $(W_{j})_{j \in J}$ be the open covering of $M$
from above and recall the  trivializations
\[
   \psi_{i,j}: F_i|_{W_j} \cong D^n \times q_i\K \, .
\]
For each $i$ and $j$, the connection $\nabla_i$ induces a smooth connection $1$-form
in $\Gamma(T^* D^n \otimes \End_{\K} (F_i)) \iso \Omega^1(D^n;q_i\K
q_i)$ using the trivialization
$\psi_{i,j}$. Since the connection is a connection of Hilbert
$\K$-modules, the values consist actually of skew-adjoint elements of
$q_iK q_i$. Using the canonical chart on $D^n$, we consider these
connections as smooth functions
\[
      \omega_{i,j} :  D^n \to (q_i\K q_i)^n
\]
and as such they have  $C^1$-norms which are uniformly bounded in $i$ and $j$.
This follows from Lemma \ref{estimate} and the curvature assumption
on the sequence $(\nabla_i)$.  In particular, the functions $\omega_{i,j}$
can be assembled to Lipschitz continuous functions
\[
      \omega^L_{j} : D^n \to (qAq)^n \, .
\]
The above isomorphism
$L \cong V$ gives rise to induced connection forms
\[
      \omega^V_{j}: D^n \to (qAq)^n
\]
equal to
\begin{equation*}
     ((\omega^V_j)(x))_{\nu} = \tau_j(x) \circ ((\omega^L_j)(x))_{\nu}
     \circ \tau_j(x)^*
   \end{equation*}
with $\nu = 1, \ldots , n$. Unfortunately, the functions $\omega^V_{j}$ need
not be smooth.

\medskip

We choose $\epsilon > 0$ so small that the $\epsilon$-neighborhood of $D^n$
in $\R^n$ is still mapped diffeomorphically to an open subset of $M$ by $\phi_j$. Now,
define a smooth function $\widetilde{\omega}^V_j : D^n \to (qAq)^n$ by
\begin{equation*}\label{eq:mollifying}
   \widetilde{\omega}^V_j (x):= \int_{D^n} \delta_{\epsilon}(x-t) \omega^V_j(t) dt
 \end{equation*}
using the Bochner integral and a smooth nonnegative bump function
$\delta_{\epsilon}:D^n \to \R$ of total integral $1$ whose support is contained
in the $\epsilon$-ball around $0$. We consider $\omega^V_j$ as
a smooth connection $1$-form in $\Gamma(T^* D^n; \End_A(V)) = \Omega^1(D^n; qAq)$
and hence as a smooth $A$-linear connection
\[
   \nabla^{V,j} \in  \Omega^1(W_j; qAq)
\]
Let $\rho_j:M \to \R$ be
a partition of unity subordinate to the covering $(W_j)$ and
set
\[
   \nabla^V := \sum_{j \in J} \rho_j \cdot \nabla^{V,j} \, .
\]
Then $\nabla^V$ is a Hilbert $A$-module connection on $V$, since the forms
$\widetilde{\omega}_j^V$ are still skew-adjoint. Since it preserves the
$A$-module structure, it also preserves the blocks $V\cdot A_i$.
We claim that it has the asymptotic curvature properties stated in Theorem \ref{main}. We denote by
\[
     \omega^L_{i,j}: D^n \to (q_i\K q_i)^n\, , ~~ \omega^V_{i,j}: D^n
     \to (q_i\K q_i)^n\, , ~~
     \widetilde{\omega}^V_{i,j}: D^n \to (q_i\K q_i)^n
\]
the connection forms that are induced by the projection $p_i\colon A
\to A_i = \K$, i.e.~$\omega^L_{i,j}=p_i\omega^L_j p_i$, etc. By construction,
\[
     \omega^L_{i,j} = \omega_{i,j} \, .
\]
Because $L$ and $V$ are Lipschitz isomorphic with a global Lipschitz constant
on $M$ (with respect to the covering of $M$ by the subsets $ W_{j}$), there is a constant
$C$ such that we have estimates of $1$-norms
\[
    \| \omega^V_{i,j}\|_1 \leq C \cdot \| \omega^L_{i,j}\|_1 = C \cdot \| \omega_{i,j}\|_1 \, .
\]
Furthermore,
\[
   \widetilde{\omega}^V_{i,j}=  \int_{D^n} \delta_{\epsilon}(x-t) \omega^V_{i,j} (t) dt \, .
\]
The formula shows that we get pointwise bounds on
$\widetilde{\omega}^V_{i,j}$ and its derivatives up to order $d$
in terms of the
sup-norm of the fixed function $\delta_{\epsilon}$ and its derivatives
up to order $d$ and of the $L^1$-norm of $\omega^V_{i,j}$.
Since the curvature of
$\widetilde{\omega}_{i,j}^V$ is in local
coordinates given by certain derivatives up to order $1$ of
$\widetilde{\omega}^V_{i,j}$ and because the derivatives
of the functions $\rho_j$, the derivatives of the transition
functions for the bundle $V$ and the derivatives of the chart
transition functions (with respect to the cover
$M \subset \cup W_j$) are globally bounded, the claim about
the asymptotic behaviour of the connection $\nabla^V$  follows.

\begin{rem} An alternative construction of $\nabla^V$ consists  of
assembling the given connections $\nabla_i$ to a Lipschitz connection
on $L$ inducing a Lipschitz connection on $V$. This is  then smoothed
to yield the desired connection $\nabla^V$. Our argument given before
avoids the discussion of Lipschitz connections.
\end{rem}

\begin{remark}
  We have been careful to write down the rather explicit connection
  $\nabla^V$ with its curvature properties, because this is used
  in Proposition
  \ref{prop:flat_inf_dim_bundle} to show that a suitable quotient
  bundle is \emph{flat}, which is the main ingredient in the proof of
  our main result Theorem \ref{index}.

Alternatively, one could use a different argument to show directly 
  that this quotient
  bundle admits a flat connection. We thank Ulrich Bunke for pointing
  out that this could be done by studying the parallel
  transport on the path groupoid of our manifold.
\end{remark}

\section{Almost flat bundles and index theory}
\label{sec:indextheory}

This section provides a link between the construction
from the last section and the index theory for Dirac operators.

Let $M^{2n}$ be a closed oriented Riemannian manifold of even dimension and
let $S \to M$ be a complex Dirac bundle equipped with a hermitian
metric and a compatible connection (cf.~\cite[Definition 5.2]{LM}).
As usual, Clifford multiplication with the complex volume element
$i^n\omega_{\C}$ induces a splitting $S^{\pm} \to M$  into $\pm 1$
eigenspaces. The corresponding Dirac type operator
\[
   D: \Gamma(S^{+}) \to\Gamma(S^{-})
\]
has an index in $K_0(\C) \cong \Z$. Denoting the
universal cover of $M$ by $\widetilde{M}$ and using the usual
representation of $\pi_1(M)$ on $C^*_{\max} \pi_1(M)$, the
maximal real $C^*$-algebra of $\pi_1(M)$, we obtain the
flat Mishchenko-Fomenko line bundle
\[
    E := \widetilde{M} \times_{\pi_1(M)} C^*_{\max} \pi_1(M) \to M \, .
\]
The twisted Dirac type operator
\[
    D \otimes {\rm id} : S^+ \otimes E \to S^- \otimes E
\]
has an index  (cf.~\cite{Ro3})
\[
  \alpha_S(M) \in K_0(C^*_{max} \pi_1(M)) \, .
\]
In order to  detect the non-triviality of $\alpha_S(M)$ in certain cases,
we use a $C^*$-algebra morphism
\[
    C^*_{\max}(\pi_1(M)) \to Q
\]
where $Q$ is another $C^*$-algebra whose $K$-theory can be understood
explicitely and study the image of $\alpha_S(M)$ under the induced map
in $K$-homology.

First, we recall the following
universal property of $C^*_{max} \pi$ for a discrete group $\pi$:
Each involutive multiplicative map
\[
   \pi \to Q
\]
with values in the unitaries of some unital $C^*$-algebra can be extended
to a unique $C^*$-algebra morphism
\[
    C^*_{max}(\pi) \to Q \, .
\]

We now prove a naturality property of indices of twisted Dirac operators.
In the following, we always use the maximal tensor product.

\begin{lem} \label{prep} Let $M$ be a compact oriented manifold of even dimension,
$S\to M$ be a Dirac bundle, $F$ and $G$ be $C^*$-algebras and $\psi:F\to G$
be a $C^*$-algebra morphism. Further, let $X$ be a Hilbert $F$-module bundle
on $M$. We define the Hilbert $G$-module bundle
\[
        Y := X \otimes_{\psi} G \, .
\]
Let $[D_X] \in K_0(F)$ and $[D_Y] \in K_0(G)$ be the indices of
the twisted Dirac type 
operators
\begin{eqnarray*}
   D_X:  \Gamma(S^+ \otimes X) & \to & \Gamma(S^- \otimes X) \\
   D_Y:  \Gamma(S^+ \otimes Y) & \to & \Gamma(S^- \otimes Y)  \, .
\end{eqnarray*}
(with an arbitrary $F$-module connection on $X$ and $G$-module
connection on $Y$). Then we have
\[
     [D_Y] = \psi_* ([D_X]) \, .
\]
\end{lem}

\begin{proof} This follows from functoriality of Kasparov's $KK$-machinery.
We denote by $[D] \in KK(\C(M), \C)$ the $KK$-element
defined by the Dirac operator $D: \Gamma(S^+) \to \Gamma(S^-)$ and by
\[
   [X] \in KK(\C, \C(M) \otimes F)
\]
the $KK$-element represented by the Kasparov triple $(\Gamma(X), \mu_X , 0)$,
where $\mu_X: \C(M) \otimes F \to \B(\Gamma(X))$ is the map induced
by the right $F$-module structure on $X$. Using the Kasparov intersection
product, we  get
\[
  [D_X] =  [X] \otimes_{\C(M)} [D] \in KK(\C ,F)  \\
\]
and  $[D_Y] \in KK(\C,G)$ is equal to
\[
  [(\Gamma(X \otimes_{\psi} G), \mu_Y , 0)] \otimes_{\C(M)} [D]  =
  [(\Gamma(X) \otimes_{\psi} G, \mu_Y , 0)] \otimes_{\C(M)} [D] \, .
\]
By definition,
\[
    \psi_*[(\Gamma(X) , \mu_F, 0)] = [\Gamma(X) \otimes_{\psi} G, \mu_Y, 0]
\]
and our claim follows from the naturality of the Kasparov intersection
product. For more details on the connection between the $KK$-description
of the index and the usual definition in terms of kernel and cokernel,
compare e.g.~\cite{Sch}.
\end{proof}

\begin{remark}\label{rem:K1index}
  Of course, in the situation of Lemma \ref{prep} there is also a
  corresponding index for odd dimensional
  manifolds, taking values in $K_1(F)$, with the corresponding
  properties. 
\end{remark}

\begin{cor} Let
\[
    \pi_1(M) \to \U(d)
\]
be a finite dimensional and unitary representation
with induced $C^*$-morphism  $\psi: C^*_{\max} \pi_1(M)  \to \Mat(\C,d)$. Let
\[
  \psi_*: K_0(C^*_{\max} \pi_1(M)) \to K_0(\Mat(\C,d)) \cong K_0(\C) = \Z
\]
be the map induced by $\psi$. Then $\psi_*(\alpha_S(M))$ coincides with the
index of the Dirac type operator $D$ twisted by the bundle
\[
     \widetilde{M} \times_{\pi_1(M)} \C^d \to M \, .
\]
\end{cor}
Here we used the following well known instance of Morita equivalence:
The index of $D$ twisted with the Hilbert $\C$-module bundle
(i.e.~vector bundle)
$\widetilde{M}\times_{\pi_1(M)} \C^d$ is equal to the
$\Mat(\C,d)$-index of $D$ twisted with the Hilbert $\Mat(\C,d)$-module
bundle $\widetilde{M}\times_{\pi_1(M)} \Mat(\C,d)$.

Unfortunately, because the higher Chern classes of finite dimensional flat
bundles vanish (using Chern-Weil theory), the element $\psi(\alpha_S(M))$
is simply equal to $d\cdot \ind(D)$.

We will now use the construction of Section \ref{sec:assembling} in order
to get a useful {\em infinite dimensional}
holonomy representation of $\pi_1(M)$.

Let $(P_i, \nabla_i)_{i \in \N}$ be a sequence of almost
flat vector bundles on $M$  and let $V \to M$ be the smooth
Hilbert $A$-module bundle constructed in Theorem \ref{main}.
Let
\[
    A' = \overline{\bigoplus_{i=1}^{\infty} \K} \subset  A
\]
be the closed two sided ideal consisting of sequences of
elements in $\K$ that converge to $0$ and let
\[
     Q := A / A'
\]
be the quotient $C^*$-algebra. The bundle
\[
    W:= V / (V \cdot A') \to M
\]
is a smooth Hilbert $Q$-module bundle with fiber $\overline{q} Q$.  Here
$\overline{q}$ is the image of the projection $q\in A$ in $Q$. The connection
$\nabla^V$ induces a connection on $W$. The
following fact follows immediately from the construction
of the bundle $V$.

\begin{prop}\label{prop:flat_inf_dim_bundle} The curvature form
\[
    \Omega_V \in \Gamma(\Lambda^2 M \otimes \End_{A}(V))
\]
can be considered as a form in
\[
    \Gamma(\Lambda^2(M) \otimes \Hom_{A}(V,V\cdot A')) \, .
\]
As a consequence, the induced connection on $W$ is flat.
\end{prop}

  Fixing a base point $x \in M$ and an isomorphism of the fiber
  $W_x\iso \overline{q}Q$, the holonomy around loops based at
  $x$ gives rise to a multiplicative involutive map
\[
\pi_1(M,x) \to \Hom_{Q}(W_x ,W_x)\iso\Hom_Q(\overline{q}Q,\overline{q}Q) = \overline{q}Q\overline{q}
\]
with values in the unitaries of the subalgebra
$\overline{q}Q\overline{q}$ of $Q$ and hence by composition to a map of
$C^*$-algebras
\[
\phi_1\colon C^*_{\max} \pi_1(M) \to \overline{q}Q\overline{q} \to Q.
\]
Let $\phi_2\colon C^*_{max}\pi_1(M)\to Q$ be the homomorphism obtained
by the same construction, but now applied to the sequence
$P_i':=M\times U(d_i)$ of trivial bundles, with trivial (and hence
flat) connections. Define
\begin{equation*}
  \phi_*:=(\phi_1)_*-(\phi_2)_*\colon K_0(C^*_{max}\pi_1(M))\to K_0(Q). \label{eq:phi}
\end{equation*}

It is not difficult to compute the $K$-theory of $A$ and $Q$.
Throughout the following argument, we work with the usual fixed isomorphism
\[
    \Z = K_0(\K) \, .
\]
Recall also that
$K_1(\K)=0$. Since $K$-theory commutes with  direct limits, we
obtain an isomorphism $K_0(A')\iso \oplus_{i=1}^\infty \Z$
and $K_1(A')=0$ (recall that $A'$ is the ideal $\oplus_{i=1}^\infty
\K$ in $A$).

\begin{prop} \label{calc} Let
\[
     J \subset \prod_{i \in \N}  \Z
\]
be the subgroup consisting of  sequences
with only finitely many nonzero elements,
i.e.~$J=\oplus_{i=1}^\infty \Z$.
Then we have
\begin{eqnarray*}
  K_0(A) & \cong & \prod_{i \in \N}  \Z\, , \\
  K_0(Q) & \cong & ( \prod   \Z )/ J\, .
\end{eqnarray*}
Under the above isomorphisms, the natural map
$K_0(A) \to K_0(Q)$ corresponds to the projection
\[
   \prod  \Z  \to ( \prod    \Z) / J \, .
\]

\end{prop}

\begin{proof} Observe that the projections
  $\{(p_{n_i})_{i\in\naturals} \mid n_i\in\naturals\}$ form an
  (uncountable) approximate unit of $A$, where
  $p_n\in\K$ is the standard projection of rank $n$. Consequently, $A$
  is stably unital in the sense of \cite[Definition 5.5.4]{B}. By
  \cite[Proposition 5.5.5]{B}, elements in $K_0(A)$ are represented by formal
differences of projections in
\[
  \Mat_\infty(A) = \Mat_\infty\left( \prod \K\right) \, ,
\]
where $\Mat_\infty$ is the union of all the $\Mat_r$. The main point
of this stable unitality is that we don't have to adjoin a unit to $A$ in
order to compute $K_0$. By projecting
to the different ``coefficients'' $A_i$ we get an
induced map
\[
   \chi:K_0(A) \to \prod_{i=1}^{\infty} K_0(\K) = \prod  \Z  \, ,
\]
Writing down appropriate projections, we see that  $\chi$ is
surjective. For the injectivity
of $\chi$, consider two projections $P,Q \in \Mat_r(A)$ such
that for all $i\in \N$ the components $P_i, Q_i \in \Mat_r(\K)$ are
equivalent, where the subscript $i$ indicates application of the
projection $A = \prod_i \K \to A_i=\K$ onto the $i$th factor. We
get a family of partial isometries $V_i \in \Mat_{r}(\K)$
such that
\[
    V_iV_i^*= P_i\qquad V_i^* V_i = Q_i \, .
    \]
Because all the matrices $V_i$ have norm $1$, they can be assembled
to a partial isometry $V \in \Mat_{r}(A)$ such that
$VV^*=P$ and $V^*V = Q$.

\medskip

The calculation of $K_0(Q)$ uses the exact sequence
\[
   K_0(A') \stackrel{\iota_*}{\ra} K_0(A) \stackrel{\pi_*}{\ra} K_0(Q)
   \to K_1(A')
   \]
  induced by the short exact sequence
  \begin{equation*}
0\to A'\xrightarrow{\iota}
  A\xrightarrow{\pi}Q\to 0,
\end{equation*}
where $\iota:A' \to A$ and $\pi: A \to Q$ are the obvious maps.
Since $K_1(A')=0$,  $\pi_*$ is surjective.

The inclusion of the first $k$ summands
\[
   \K \oplus \ldots \oplus \K \to A
\]
induces an injective map
\[
   K_0( \K \oplus \ldots \oplus \K ) \to K_0(A)
\]
that can be identified with the inclusion
\[
   \Z \oplus \ldots \oplus \Z \to \prod  \Z
\]
onto the first $k$ factors. The map $\iota_*$ is now given by passing
to the inductive limit of the last map, and this finishes
our calculation of $K_0(Q)$.
\end{proof}

\begin{rem} In a similar way, it can be shown that
\[
  K_1(A) = K_1(Q) = 0 \, .
\]
\end{rem}

Now we can formulate the following important fact which shows
that the asymptotic index theoretic information of the sequence of
almost flat bundles
$(P_i)$ is completely contained in $\alpha_{S}(M)$.

\begin{thm} \label{index} For all $i \in \Z$, define
  \begin{equation*}
z_i := \ind(D_{E_i}) - d_i\ind(D)\in K_0(\C) = \Z,
\end{equation*}
the index of the Dirac type operator $D$ twisted by the virtual bundle
$E_i- \underline{\C}^{d_i}$
where $E_i \to M$ is the $d_i$-dimensional unitary vector bundle
with connection induced by the connection $\nabla_i$ on $P_i$,
$S^{\pm} \otimes E_i$ is  equipped with the product connection and
$\underline{\C}$ is the
trivial bundle. Then the element
\[
      \phi_*(\alpha_S(M)) \in K_0(Q)
\]
is represented by
\[
    (z_i) \in \prod \Z  \, .
\]
\end{thm}

\begin{proof} The idea of the proof is to study the image of a
(computable) index of the Dirac operator twisted
with a {\em non-flat} bundle of $A$-modules over $M$ under the canonical map
\[
   K_0(A) \to K_0(Q) \, .
\]
Using Lemma \ref{prep}, this index turns out to be  equal to the
index of $D$ twisted with a {\em flat} bundle which is induced by the
given holonomy representation of $\pi_1$ on $Q$. Therefore it is
equal to $\phi_*(\alpha_S(M))$.

In order to make this idea precise, we consider the bundle of
$A$-Hilbert modules $V \to M$ constructed in
Theorem \ref{main} and the element
\[
   [D_V] \in KK(\C,A)\cong K_0(A) =  \prod\Z
\]
represented by the Dirac operator
\[
   D_V: \Gamma(S^+ \otimes V) \to \Gamma(S^- \otimes V)
\]
on $M$. For $i \in \N$ let
\[
   p:A \to \K
\]
be the projection onto the $i$th factor. By Lemma \ref{prep}, the induced map
\[
    p_*:K_0(A) \to K_0(\K) \cong  \Z
\]
sends $[D_V]$ to the index of $D_{P_i \times_{\U(d_i)} \K}$. Hence,
\[
    p_*([D_V]) = \ind(D_{E_i}).
\]
If we carry out the same construction with the trivial $U(d_i)$-bundle
$P_i'$ we obtain
\begin{equation*}
p_*([D_{V'}]) = d_i\ind(D),
\end{equation*}
and it follows that
\[
   [D_V] - [D_{V'}] = (z_1,z_2,z_3, \ldots) \, .
\]
Under the canonical map
\[
   K_0(A) \to K_0(Q) \, ,
\]
the element $[D_V]$ is mapped to the element represented by the index
of the Dirac operator $D$ twisted with the flat bundle
\[
    W = \widetilde{M} \times_{\pi_1(M)} Q
\]
using the holonomy representation $\phi_1$ constructed from the $(P_i)$ of
$\pi_1(M)$ on $Q$.
This element coincides with $(\phi_1)_*(\alpha_S(M))$.
In a similar way, $(\phi_2)_*(\alpha_S(M))$ is the image of
$[D_{V'}]$ under the canonical map $K_0(A)\to K_0(Q)$ and it
remains to take the difference in order to finish the
proof of Theorem  \ref{index}.
\end{proof}

The reason for using the virtual bundles $E_i - \C^{d_i}$ will
become apparent in the applications described in the next sections.

\section{Enlargeability and universal index}
\label{sec:enlargeability}

For a closed spin manifold $M^{2n}$ of even dimension, we consider the Dirac
bundle $S \to M$ given by the complex spinor bundle on $M$.
In this case, we define
\[
   \alpha_{max}(M) \in K_{2n}(C^*_{max} \pi_1(M)) = K_0(C^*_{max} \pi_1(M))
\]
to be equal to $\alpha_S(M) \in K_0(C^*_{max} \pi_1(M))$ (cf.~Section
\ref{sec:indextheory}). If the dimension of $M$ is odd, note that
\begin{equation*}
  K_0(C^*_{max}(\pi_1(M)\times \integers)) =
  K_0(C^*_{max}\pi_1(M))\tensor 1
  \oplus K_1(C^*_{max}\pi_1(M))\tensor e,
\end{equation*}
using the exterior Kasparov product 
\begin{equation*}
K_*(C^*_{max}\pi_1(M))\tensor
K_*(C^*\integers)\to K_*(C^*_{max}(\pi_1(M)\times\integers))
\end{equation*}
with the canonical
generators $1\in K_0(C^*\integers)$ and $e\in
K_1(C^*\integers)$.
 Using this splitting, we  define
$\alpha_{max}(M)\in K_1(C^*_{max}\pi_1(M))$ by requiring that
\begin{equation*}
\alpha_{max}(M)\tensor e = \alpha_{max}(M\times S^1).
\end{equation*}
This is consistent with the direct definition of $\alpha_{max}(M)$ alluded
to in Remark \ref{rem:K1index} and the product formula \cite[Theorem
9.20]{Stolz(1998)}
\begin{equation*}
  \alpha_{max}(M\times S^1) = \alpha_{max}(M)\tensor e,\text{ with }e=\alpha_{max}(S^1).
\end{equation*}

The following fact is well known and can be proven in the usual
way by an appropriate Weitzenb{\"o}ck formula.

\begin{prop} Let $M$ be a closed spin manifold. If $M$ admits a metric of
positive scalar curvature, then
\[
   \alpha_{max}(M) = 0 \, .
\]
\end{prop}

\begin{thm} \label{enlargeable} Let $M^m$ be an enlargeable or
  area-enlargeable spin manifold. Then
\[
    \alpha_{max}(M) \neq 0 \in K_m(C^*_{max}(\pi_1(M))) \, .
\]
\end{thm}

\begin{proof} We first show how we can reduce to the case that $M$ has
  even dimension. If not, consider the commutative diagram
  \begin{equation*}
    \begin{CD}
      K_m(M) @>{\times[S^1]}>> K_{m+1}(M\times S^1)\\
      @VVV @VVV\\
      K_m(B\pi_1(M)) @>{\times[S^1]}>> K_{m+1}(B\pi_1(M)\times
      B\integers)\\
      @VVV @VVV\\
      K_m(C^*\pi_1(M)) @>{\times\alpha_{max}(S^1)}>> 
 K_{m+1}(C^*(\pi_1(M)\times\integers))\, .
    \end{CD}
  \end{equation*}
Here we use the fact that the $\alpha$-index is multiplicative with
respect to the exterior Kasparov product (note that  $B\integers=S^1$ and
$C_{max}^*(\pi_1(M)\times\integers)=C_{max}^*\pi_1(M)\tensor C^*\integers$),
compare \cite[Theorem 9.20]{Stolz(1998)}. Since $M\times S^1$ is (area)-enlargeable if
$M$ is, and the image of $\alpha_{max}(M)$ under the bottom horizontal arrow
is $\alpha_{max}(M\times S^1)$, it suffices to treat non-vanishing of this
invariant for even dimensional area-enlargeable spin manifolds.

Therefore, we assume that $M$ is of even dimension $2n$
so that $\alpha_{max}(M)$ can
be considered as an element in $K_0(C^*_{max} \pi_1(M))$.

Because $M$ is area-enlargeable, there is a sequence of almost
flat principal unitary bundles $(P_i)$ on $M$ such that the Chern classes
in $H^*(M; \Z)$ of the associated (finite dimensional) complex vector bundles
$E_i$ 
satisfy
\begin{eqnarray*}    
   c_{\nu}(E_i) & = &   0 \, , {\rm~if~}  0 <  \nu < n \\
   { \langle c_{n}(E_i), [M] \rangle} &\ne & 0,\text{ if } \nu=n \, .
\end{eqnarray*}
Such a sequence can be constructed as follows: Because the Chern character
\[
   K^0(S^{2n}) \otimes \Q \to H^{even}(S^{2n}; \Q)
\]
is an isomorphism, there is a vector bundle
\[
   E \to S^{2n}
\]
with
\[
   c_n(E) \neq 0 \in H^{2n}(S^{2n}; \Z)\, .
\]
Now let $i \in \N$ and choose a finite  covering $\overline{M} \to
M$ with covering group $G$ such that there is a $\frac{1}{i}$-area contracting
map $\psi: \overline{M} \to S^{2n}$ of nonzero degree. Passing 
to a finite cover of $\overline{M}$ if necessary, we can 
assume without loss of generality that the covering
$\overline{M} \to M$ is regular. The $G$-action  on
$\overline{M}$ can be extended to an action of this group
on
\[
    \bigoplus_{g \in G} g^*(\psi^*(E))
\]
by vector bundle automorphisms. Note that the norm of the
curvature of this direct sum of bundles is not larger than the norm of
the curvature of 
$\psi^*(E)$ and this is bounded by $\frac{1}{i}$ times the norm of the
curvature of $E$ by the $\frac{1}{i}$-area contractibility. (If the map was
$\frac{1}{i}$-contractible, we would get a factor $\frac{1}{i^2}$). Let $E_i \to M$ be the quotient
vector bundle. By naturality of Chern classes we have $c_{n}(E_i) \neq 0$ and $c_{\nu}(E_i)= 0$, if
$0 < \nu < n$. The last statement is true, because the canonical map
\[
   H^*(M; \Q) \to H^*(\overline{M}; \Q)
\]
is injective, transfer followed by division by $n$ giving a splitting.

By construction, the Chern character of the virtual vector bundle
$E_i-\underline{\C}^{d_i}$  is
\begin{equation*}
  \ch(E_i-\underline{\C}^{d_i}) = C  \cdot c_n(E_i) \ne 0
\end{equation*}
with some non-zero constant $C$ (dependent of $n$). In
particular, $\ch(E_i - \underline{\C}^{d_i})$ is concentrated
in degree $2n$.

The Atiyah-Singer index formula implies that the integer
valued index in $K_0(\C) \cong \Z$ of the Dirac operator
\begin{equation*}
D_{E_i-\underline{\C}^{d_i}}\colon \Gamma(S^{+} \otimes
(E_i-\underline{\C}^{d_i})) \to  \Gamma(S^{-} \otimes
(E_i-\underline{\complexs}^{d_i}))
\end{equation*}
is equal to
\begin{equation*}
 \innerprod{\hat{\mathcal{A}}(TM) \cup
    \ch(E_i-\underline{\C}^{d_i}),[M]} = C \cdot
    \innerprod{c_n(E_i),[M]} 
\ne 0
\end{equation*}
where $\hat{\mathcal{A}}$ denotes the total $\hat{A}$-class. Now, Theorem \ref{index} implies our assertion.
\end{proof}

\section{On a question by Burghelea}

\noindent {\bf Question} (\cite[Problem 11.1]{Schu})  ``If $M^n$ is an
enlargeable manifold and
\[
   f : M \to B\pi_1(M)
\]
induces an isomorphism on the fundamental groups, does $f_*$
map the fundamental class of $H_n(M;\Q)$ non-trivially? Is the
converse statement true?''

The next theorem answers the first question affirmatively. We can even
drop any spin assumption on $M$ or its universal cover. At the 
end of this section, we will show by an example that the converse
of Burghelea's question in its stated form must be  answered in the negative.  

\begin{thm} \label{Burg1} Let $M$ be an enlargeable or
  area-enlargeable manifold of
  dimension $m$. Then
\[
    f_*([M]) \neq 0 \in H_{m}(B \pi_1(M); \Q) \, .
\]
\end{thm}

We can assume that $M$ is connected. We first reduce to the case that
$M$ has even dimension $2n$. Else, observe that $M\times S^1$ also is
enlargeable and we have the commutative diagram
\begin{equation*}
  \begin{CD}
    H_m(M;\rationals) @>{\times [S^1]}>> H_{m+1}(M\times
    S^1;\rationals)\\
    @VVV @VVV\\
    H_m(B\pi_1(M);\rationals) @>{\times [S^1]}>> 
H_{m+1}(B(\pi_1(M)\times\integers);\rationals)\, ,
\end{CD}
\end{equation*}
where the image of the
fundamental class of $M$ is mapped to the image of the fundamental
class of $M\times S^1$ under the bottom horizontal map.

Given $M$ of dimension $2n$, we choose a sequence of almost flat bundles $(E_i)$ as
in the proof of Theorem \ref{enlargeable}.
Now consider the commutative diagram
\begin{equation}\label{eq:f_and_beta}
\begin{CD}
K_0(M)\otimes \Q @> f_*>> K_0(B \pi_1(M)) \otimes \Q   @>\beta >> K_0(Q) \otimes \Q \\
@V \ch V V       @V \ch VV                            @V=VV   \\
H_{\rm even}(M;\Q) @>f_*>> H_{\rm even}(B\pi_1(M);\Q) @>>> K_0(Q) \otimes \Q
\end{CD} \, 
\end{equation}
In this diagram, the map $\beta$ is induced by the  composition of the
assembly map
\[
   \mu: K_0(B \pi_1(M)) \to K_0(C^*_{max} \pi_1(M))
\]
with the map $K_0(C^*_{max} \pi_1(M)) \to K_0(Q)$ which is induced by the
map $\phi$ defined in \eqref{eq:phi} and
the almost flat sequence $(E_i)$.  Furthermore, $\ch$
denotes the homological Chern character.

We need the following description of
the $K$-homology $K_0(M)$ from \cite{Kes}, Definition 2.6 and Lemma 2.8.

\begin{prop} Let $X$ be a connected $CW$-complex.
Elements in $K_0(X)$ are represented by triples $(N,S,u)$, where $N$ is a closed oriented Riemannian manifold of even dimension
(consisting of components of possibly different dimension), $S \to M$
is a Dirac bundle on $M$ and $u: N \to X$ is a continuous map.
Two such triples are identified, if they are equivalent under the equivalence
relation  generated by
direct sum/disjoint union, bordism and vector bundle modification.

\end{prop}

Using vector bundle modification and the bordism relation introduced above, we can assume that in the
triple $(N, S, u)$ above, the manifold $N$ is connected.

Now let $M$ be the given manifold and let $(N,S,u)$ represent
an element $k \in K_0(M)$. Let $D:\Gamma(S^+) \to \Gamma(S^-)$ be
the Dirac type operator associated to the Dirac bundle $S$.

\begin{lem} \label{lemma1} The element $\beta \circ f_*(k)$ (compare
  \eqref{eq:f_and_beta}) is represented
by
\[
    (z_1, z_2, \ldots) \in \prod \Z = K_0(A)
\]
where
\[
   z_i = \ind(D_{u^*(E_i) - \underline{\C}^{d_i}})
\]
is the index of $D$ twisted by the virtual bundle $u^*(E_i) -
\underline{\C}^{d_i}$.
\end{lem}

The proof of this statement is analogous to the proof of
Theorem \ref{enlargeable}. Note that the Kasparov $KK$-theory
element in $K_0(M)$ represented by $(N,S,u)$ is equal to
\[
    u_*([D])
\]
where $[D] \in KK(\C(N), \C)$ is the $KK$-element induced by $D$
(cf.~the explanations before Example 2.9.~n \cite{Kes}).
In particular, the element
\[
    \mu \circ f_*(k) \in K_0(C^*_{max} \pi_1(M))
\]
is given as the index of the Dirac operator $D$ twisted by $u^*(E)$,
where $E \to M$ is the Mishchenko-Fomenko line bundle on $M$
with fibre $C^*_{max} \pi_1(M)$.

\begin{lem} \label{lemma2} The element $\ch(k) \in H_{even}(M; \Q)$ is equal to
\[
  (-1)^n \cdot u_*\Big( \big( p_{!}  \ch(\sigma(D) ) \cup \Todd(TN \otimes \C)
\big) \cap [N] \Big) \, ,
\]
where $\sigma(D)$ is the $K$-theoretic symbol class, $\Todd$ is
the total Todd class and
\[
   p_!: H_c^*(TN;\Q) \to H^{*-\dim(N)}(N;\Q)
\]
is the Gysin map induced by the canonical projection $p : TN \to N$.
\end{lem}

\begin{proof} In a first step, one shows that the assignment
\[
  \omega:  (N, S , u) \mapsto (-1)^n \cdot u_*\Big( \big( p_! \ch(\sigma(D)) \cup \Todd (TN \otimes \C)\big) \cap [N] \Big)
\]
is compatible with the equivalence relation imposed on the
triples $(N, S, u)$ used in the definition of $K_0(M)$. For
disjoint union and bordism, this is straightforward. The
invariance under vector bundle modification uses the same
calculation as in Section 7 of \cite{EG}, p.~64. Consequently, $\omega$
induces an additive map
\[
   K_0(M) \otimes \Q \to H_{even}(M; \Q) \, .
\]
In order to prove that this map is indeed equal to the homological
Chern character, it is enough to consider triples $(N, S, u)$,
where $N$ is a $\Spc$-manifold and $S$ is the canonical spinor
bundle on $N$ (cf.~\cite[2.3]{Kes}). But in this special case
an explicit calculation shows that
\[
  (-1)^n \cdot p_! \ch(\sigma(D)\big) \cup \Todd(TN \otimes \C) = e^{\frac{1}{2} c}
\cdot \hat{\mathcal{A}}(TN)
\]
where $c \in H^2(N;\Q)$ is the first Chern class of the complex
line bundle associated to the $\Spc$-structure on $N$. Now one
uses the calculation of the homological Chern character in \cite[4.2]{J}.
\end{proof}

We continue the proof of Theorem \ref{Burg1}. Let $(N,S,u)$ be
a triple (with connected $N$) representing an element in
$K_0(M) \otimes \Q$ which under the homological Chern
character is mapped to $q \cdot [M] \in H_{2n}(M ;\Q)$.
Here, $q$ denotes an appropriate nonzero rational number. As
before, let $D:\Gamma(S^+)
\to \Gamma(S^-)$ be the associated Dirac type operator.
We will show that
\[
    \beta \circ f_*([N,S,u]) \neq 0 \in K_0(Q)
\]
which implies the assertion of Theorem \ref{Burg1} by the commutativity
of diagram \eqref{eq:f_and_beta}.

Using Lemma \ref{lemma1},
\[
  \beta \circ f_*([N,S,u]) \in \prod \Z \big/ \bigoplus \Z = K_0(Q)
\]
is represented
by the sequence $(z_1, z_2, \ldots) \in \prod \Z$ where
\[
   z_i = \ind(D_{f^*(E_i) - \underline{\C}^{d_i}}) \, .
\]
By the Atiyah-Singer index theorem, this index is given by the zero
dimensional component of the homology class
\[
 (-1)^n \Big( p_!\ch(\sigma(D)) \cup \Todd(TN \otimes \C) \cup \ch(u^*(E_i) - \underline{\C}^{d_i})\Big )  \cap [N]  \in H_*(N; \Q)\, .
\]
On the other hand, because
\[
    u_*:H_*(N;\Q) \to H_*(M;\Q)
\]
induces an isomorphism in degree $0$, the number $z_i$ is equal to  the
zero dimensional component of
\begin{eqnarray*}
\lefteqn{u_* \bigg(  (-1)^n \Big( p_!\ch(\sigma(D)) \cup \Todd(TN \otimes \C) \cup \ch\big(u^*(E_i) - \underline{\C}^{d_i}\big)\Big )\cap [N] \bigg) = }
\\
 & &   \ch(E_i-\underline{\C}^{d_i}) \cap u_* \Bigg( (-1)^n \Big( p_!\ch(\sigma(D)) \cup \Todd(TN \otimes \C) \Big) \cap [N] \Bigg) = \\
 & & q \cdot \big( \ch(E_i - \underline{\C}^{d_i}) \cap [M] \big) \, .
\end{eqnarray*}
The last equality uses Lemma \ref{lemma2}.
Hence,
\[
   z_i = q \cdot C \cdot \innerprod{ c_n(E_i) , [M] } \in q \cdot (\Z \setminus \{0\})
\]
by the construction of the sequence $E_i$ (the constant $C$ was
introduced at the end of Section  \ref{sec:enlargeability}).
It follows that $\beta \circ f_*([N,S,u]) \neq 0$ and
the proof of Theorem \ref{Burg1} is complete.

The following Lemma prepares the construction of 
an example showing that the converse of Burghelea's question 
must be answered in the negative. 

\begin{lem} \label{converse} For every natural number 
$n > 0$, there is a finitely presented group $G$ without
proper subgroups of finite index and such that 
\[ 
   H_n(G;\Q) \neq 0 \, .
\]
\end{lem}

\begin{proof} The proof is modeled on a similar construction 
in \cite{BDH} (cf.~Theorem 6.1 and the following remarks in
this reference). Let 
\[
  K_1:=\langle a,b,c,d~|~a^{-1}ba = b^2, b^{-1}cb= c^2,
  c^{-1}dc=d^2, d^{-1}ad= a^2 \rangle
\]
be the Higman group \cite{Hig}. This is a finitely presented 
infinite group without nontrivial finite quotients (and hence 
without proper subgroups of finite index). By 
\cite{BDH}, $K_1$ is acyclic and in particular   
\[
   \widetilde{H}_*(K_1;\Q)= 0 \, .
\]
There is an element $z \in K_1$ generating 
a subgroup $G_1 < K_1$ of infinite order. 
The amalgamated product
\[
   G_2:=K_1 *_{G_1} K_1 
\]
is finitely presented and does still have no nontrivial  
finite quotients as one checks directly with help of 
the universal property of push-outs.  
The Mayer-Vietoris sequence shows that
\[
   H_*(G_2;\Q) \cong H_*(S^2;\Q) \, .
\]
By \cite[Theorem 6.1]{BDH} the group $G_2$ embeds in the
rationally acyclic group without nontrivial finite quotients
\[
   K_2 := (K_1 \times G_1) *_{G_1} K_1 \, .
\]
Here we identify $G_1$ on the left with the second factor of $K_1 \times
G_1$. We set
\[
   G_3 := K_2 *_{G_2} K_2 \, .
\]
Again using the Meyer-Vietoris sequence,
\[
   H_*(G_3;\Q) \cong H_*(S^3;\Q) \, .
\]
This process is now carried out inductively by embedding
$G_{i}$ in the rationally acyclic finitely presented group without nontrivial
finite quotients
\[
   K_i := (K_{i-1} \times G_{i-1}) *_{G_{i-1}} K_{i-1} 
\]
and defining
\[
   G_{i+1} := K_i *_{G_i} K_i \, .
\]
The group $G:= G_n$ has then the desired properties. 
\end{proof}
        
The following theorem provides a negative answer to the 
converse of the question by Burghelea.

\begin{thm} Let $n \geq 4$ be a natural number. Then there exists 
a closed $n$-dimensional spin manifold $M$ which is not 
area-enlargeable, but whose classifying map $M \to B \pi_1(M)$ 
sends the fundamental class of $M$ to a nonzero class
in $H_n(B \pi_1(M);\Q)$. 
\end{thm}

\begin{proof} Let $G$ be the group constructed in Lemma
\ref{converse} for the number $n$. The Atiyah-Hirzebruch
spectral sequence computing the rational spin bordism 
of $BG$ (or any other space) collapses at $E_2$, and hence there 
is a closed $n$-dimensional 
spin manifold $N$ together with a map $f:N \to BG$ such that 
\[
   f_*([N]) \neq 0 \in H_n(BG,\Q) \, .
\]
Because $n \geq 4$ and because $G$ is finitely presented,
there is also an $n$-dimensional closed spin manifold
$A$ with fundamental group $G$. Let $g:A  \to BG$ denote the
classifying map and consider the map 
\[
   N \sharp A   \stackrel{f \sharp g}{\longrightarrow}
 BG \, . 
\]
This map is surjective on $\pi_1$ and sends the fundamental 
class of $N \sharp  A$ to a nontrivial class in $H_n(BG;\Q)$
(if this is not the case, simply take the connected sum with 
 more copies of $A$). Carrying out spin surgery on $N
\sharp A$ over $BG$ in order to kill the kernel of $\pi_1(f \sharp g)$, 
we obtain a spin manifold $M$ with fundamental group $G$ and
such that the classifying map $M \to BG$ sends the fundamental 
class of $M$ to a nontrivial class in $H_n(BG;\Q)$. However,
$G$ does not have any proper subgroups of 
finite index and therefore $M$ does not have any nontrivial finite covers
whatsoever. Consequently, $M$  is 
not area-enlargeable. 
\end{proof}

\begin{remark}
  It is not clear if the converse of Burghelea's question has an
  affirmative answer when working with a notion of enlargeability
  allowing infinite covers. We will address this question, and the
  relation between the corresponding Gromov-Lawson obstruction to positive
  scalar curvature \cite{GLIHES} and $\alpha_{max}$ at another
  place.
\end{remark}

\section{Concluding remarks}\label{sec:concluding-remarks}

In \cite{GL}, it is shown that enlargeability is
a homotopy invariant and is preserved under
some natural geometric constructions such as taking connected
sums or taking the cartesian product of two enlargeable
manifolds. If all manifolds under consideration are spin, then
using 
the universal property of $C^*_{max}$, one can
show by purely formal arguments that the manifolds resulting
from these constructions have nonvanishing $\alpha_{max}$.
On the other
hand, this reasoning can be used to prove the nonvanishing
of $\alpha_{max}$ in some cases that do not
seem to be accessible to the classical geometric
arguments by Gromov and Lawson. For example, using the methods
developed in this paper, one can show

\begin{prop} \label{last} Let $F$ and $M$ be connected
enlargeable spin manifolds
of even dimension. Let
\[
     F \hookrightarrow E \to M
\]
be a smooth fibre bundle admitting a spin structure and inducing a
split short exact sequence
\[
   1 \to \pi_1(F) \to \pi_1(E) \to \pi_1(M) \to 1
\]
which is equivalent to the canonical split sequence
\[
   1 \to \pi_1(F) \to \pi_1(F) \times \pi_1(M) \to \pi_1(M) \to 1 \, .
\]
Then $\alpha_{max}(E) \neq 0$ and in particular $E$ does not admit
a metric of positive scalar curvature.
\end{prop}

\bigskip\noindent
{\small Bernhard Hanke\\ Universit{\"a}t M{\"u}nchen\\
  Germany\\
  \href{http://www.mathematik.uni-muenchen.de/personen/hanke.html}{http://www.mathematik.uni-muenchen.de/personen/hanke.html}}

\medskip\noindent
{\small Thomas Schick\\ Georg-August-Universit{\"a}t G{\"o}ttingen\\
  Germany\\ \href{http://www.uni-math.gwdg.de/schick}{http://www.uni-math.gwdg.de/schick}}

\end{document}